\documentclass[12pt]{amsart}
\usepackage{latexsym,amssymb}
\usepackage{amsmath}
\usepackage{amscd}

\input xypic
\newcommand{\rar}{\rightarrow}
\newcommand{\lar}{\longrightarrow}

\newcommand{\llar}{-\kern-5pt-\kern-5pt\longrightarrow}

\newtheorem{Theorem}{Theorem}[section]

\newtheorem{proposition}[Theorem]{Proposition}

\newtheorem{Example}[Theorem]{Example}

\newtheorem{Definition}[Theorem]{Definition}
\newtheorem{question}[Theorem]{Question}

\def\ann{\mbox{\rm ann}}

\def\demo{\noindent{\bf Proof. }}
\def\depth{\mbox{\rm depth }}

\def\hdeg{\mbox{\rm hdeg}}
\def\reg{\mbox{\rm reg}}
\def\Ext{\mbox{\rm Ext}}

\def\Hom{\mbox{\rm Hom}}
\def\ker{\mbox{\rm ker}}

\def\QED{\hfill$\Box$}

\def\rank{\mbox{\rm rank}}

\def\Deg{\mbox{\rm Deg}}
\def\hdeg{\mbox{\rm hdeg}}

\def\Tor{\mbox{\rm Tor}}

\def\FF{{\bf F}}

\def\xx{{\bf x}}

\def\ff{{\bf f}}

\def\g2{{\bf g}}
\def\hh{{\bf h}}
\def\xx{{\bf x}}

\def\bbq{{\mathbb Q}}
\begin{document}

\title{{\sc Length Complexity  of Tensor Products}}

\author{Wolmer V. Vasconcelos}

\thanks{Partially supported by the NSF}
\address{Department of Mathematics, Rutgers University}
\address{110 Frelinghuysen Rd, Piscataway, NJ 08854-8019, U.S.A.}
\email{vasconce@math.rutgers.edu}
\thanks{
\endgraf
}
\date{\today}

\begin{abstract} In this paper we
 introduce techniques  to gauge the torsion of the tensor product
$A\otimes_RB$ of two finitely generated modules over a Noetherian
ring $R$. The 
 outlook   is very
different from the study of the rigidity of Tor carried out in the
   work of Auslander (\cite{Aus61}) and other authors. Here the
   emphasis in on 
   the search for  bounds for
the torsion part of $A\otimes_R B$ in terms of global invariants of $A$ and
of $B$ in special classes of modules: vector bundles and modules of
  dimension at most three.

\end{abstract}
\maketitle

\tableofcontents

\section{Introduction}

\noindent
Let  $R$ be an integral domain and $A$ and $B$ finitely generated, torsionfree
$R$-modules.   It is  a challenging task to determine whether the tensor
product 
$A\otimes_R B$ is also torsionfree, and if not, what are the
invariants of its torsion
submodule. It was a remarkable discovery by M. Auslander
(\cite{Aus61}; see also \cite{Licht66}) that for regular local rings
the torsionfreeness of $A\otimes_R B$ makes great but precise demands on the
whole homology of $A$ and of $B$. These results have been taken up by
other authors, with an important development being \cite{HW94} with
its 
extensions to various types of complete intersections.

\medskip

The focus of the study of the torsion of $A\otimes_RB$ carried out here is very
different from that of the rigidity of Tor carried out in the works
mentioned above. The overall goal is that of determining bounds for
the torsion part of $A\otimes_R B$ in terms of invariants of $A$ and
of $B$. 
Let $(R, \mathfrak{m})$ be a
Noetherian local ring, and let $A$ and $B$ be finitely generated
$R$-modules. The HomAB question, first treated in \cite{Dalili} and
further developed in \cite{DV2}, asks for estimates for the number of
generators of $\Hom_R(A,B)$ in terms of invariants of $A$ and $B$ (or
even of $R$). Since the answer in special cases depends on
cohomological properties of $A$ and $B$, it seems appropriate to
express bounds for $\nu(\Hom_R(A,B))$ in  terms of some extended
multiplicity functions $\Deg(A)$ and $\Deg(B)$ of $A$ and $B$ (see \cite{DGV}).

\medskip

Here we consider  an analogue of the HomAB
question for tensor products which we formulate as follows:

\begin{question} {\rm Can the torsion of a tensor product
 be estimated in terms of
multiplicity invariants of $A$ and $B$? In particular, how to
approach  the calculation of 
$\lambda(H^0_{\mathfrak{m}}(A\otimes_RB))$ or the related
$\nu(H^0_{\mathfrak{m}}(A\otimes_RB))$?
In this note
 we look for the
existence of 
polynomials $\ff(x,y)$ with rational coefficients depending on
invariants of $R$ such that
\[ h_0(A\otimes_R B)= \lambda(H^0_{\mathfrak{m}}(A\otimes_RB))\leq \ff(\Deg(A), \Deg(B)).\]
More generally, we look for similar  bounds for  \[
h_0(H^0_{\mathfrak{m}}(\Tor_i^R(A, B))), i\geq 0.\]
}\end{question}

We shall consider special cases of these questions but 
in  dimensions $\leq 3$, vector bundles and some classes of graded
modules. The most sought after kind of answer for the shape of
$h_0(A\otimes B)$ has the format
\[ h_0(A\otimes B)\leq c(R) \cdot \Deg(A) \cdot \Deg(B),\]
where $c(R)$ is a function depending on the dimension of $R$ and some
of its related invariants such  as its Betti numbers.
Two of such results are: (i) first 
Theorem~\ref{h0vb} that asserts that $c(R)= \dim R$ works if $R$ is a
regular local ring and $A$ is a module free on the punctured
spectrum, then (ii)
 Theorem~\ref{h0dim3} that establishes a similar result, with
$c(R)< 4$, if both $A$ and $B$ are torsionfree modules over a
$3$-dimensional regular local ring.

\medskip

The difficulties mount rapidly if both $A$ and $B$ have torsion. Thus
in Theorem~\ref{grdim2}, if $A$ is a graded torsion module over
$k[x,y]$, generated by elements of degree zero, then the best we have
achieved is to get $h_0(A\otimes A)< \hdeg(A)^6$.

\medskip
On the positive side, 
one general argument (Theorem~\ref{h0highcoh}) will show that any
polynomial bound for $h_0(A\otimes B)$ leads to a similar bound for
$h_0(\Tor_i^R(A,B))$, when the cohomological degree function
$\hdeg(\cdot)$ is used.

\section{Preliminaries}

Throughout $(R, \mathfrak{m})$ is either a Noetherian local ring or a
polynomial ring over a field. For unexplained terminology we refer to
\cite{BH}. For simplicity of notation, we often denote the tensor
product of the $R$-modules $A$ and $B$ by $A\otimes B$. For a module
$A$, $\nu(A)$ will denote its minimal number of generators. If $A$ is
a module with a composition series, we denote its length by $\lambda(A)$.
Two sources, where some of the  techniques  used with cohomological
degrees were
applied to the HomAB questions,  are \cite{Dalili} and \cite{DV2}.

\subsubsection*{Finite support}

We will use the following notation.  Set
$H^0_{\mathfrak{m}}(E)=E_0$, and $h_0(E) = \lambda(E_0)$.
More generally, if $H^i_{\mathfrak{m}}(E)$ has finite length, we set
$h_i(E)= \lambda(H^i_{\mathfrak{m}}(E))$.

\medskip

Let us begin with
the following two observations. 

\begin{proposition} \label{AmodA0} Let $(R, \mathfrak{m})$ be a Noetherian local
ring and $A, B$ finitely generated $R$-modules. Then
\begin{eqnarray*}
 h_0(A\otimes_R B)&\leq &
h_0(A)\cdot \nu(B) + h_0(B)\cdot \nu(A) + h_0(A/A_0\otimes_R
B/B_0),\\
\nu(H_{\mathfrak{m}}^0(A\otimes_R B)) &\leq & \nu(A_0)\cdot\nu(B) +
\nu(A)\cdot \nu(B_0) + \nu(H_{\mathfrak{m}}^0(A/A_0\otimes_R B/B_0). \\
\end{eqnarray*}
\end{proposition}

The other reduction involves replacing $A$ and $B$ by their direct
sum $A\oplus B$. This is allowed because the invariants of $A$ and
$B$ are essentially additive, and the analysis of $A\otimes A$
is often simpler than that of $A\otimes B$.

\subsubsection*{Big Degs}

Clearly, to estimate $h_0(A\otimes B)$, the
knowledge
$h_0(A)$ and $h_0(B)$ are far from enough. Even more strongly, data
provided from the multiplicities $\deg (A)$  and $\deg (B)$ will also fall
far short of the goal. We will focus instead on the so-called {\em
extended} or {\em cohomological} degree functions of \cite{DGV}. (See
\cite[Section 2.4]{icbook} for a discussion.)

\medskip

A part to the definition of these functions is a method to
select appropriate hyperplane sections, a process that may
require that the residue fields of the rings  be infinite. 
The notions apply to standard graded algebras or local rings, as the
following  definitions illustrate.

\begin{Definition}{\rm A {\em cohomological
degree}, or {\em extended
multiplicity function},
 is a mapping \label{Degnu} from the category $\mathcal{M}(R)$ of
 finitely generated $R$-modules, 
\[\Deg(\cdot) : {\mathcal M}(R) \mapsto {\mathbb N},\]
that satisfies the following conditions.
\begin{itemize}
\item[\rm {(i)}]  If $L = H_{\mathfrak m}^0(M)$ is the submodule of
elements of $M$ that are annihilated by a 
power of the maximal ideal and $\overline{M} = M/L$, then
\begin{eqnarray}\label{hs0}
\Deg(A) = \Deg(\overline{M}) + \lambda(L),
\end{eqnarray}
where
$ \lambda(\cdot)$ is the ordinary length function.
\item[\rm {(ii)}] (Bertini's rule) 
 If $M$ has positive depth, there is $h\in \mathfrak{m} \setminus
 \mathfrak{m}^2$,
  such that
\begin{eqnarray} \label{Bertini}
\Deg(M) \geq \Deg(M/hM).
\end{eqnarray}
\item[{\rm (iii)}] (The calibration rule) If $M$ is a Cohen-Macaulay
module, then
\begin{eqnarray} \label{calibration}
\Deg(M) = \deg(M),
\end{eqnarray}
where $\deg(M)$ is the ordinary multiplicity of  $M$.
\end{itemize}
}\end{Definition}

These functions will be referred to as big Degs. If $\dim R =0$,
$\lambda(\cdot)$ is the unique Deg function. For $\dim R=1$, $\Deg(A)=
\lambda(L) + \deg(A/L)$. When $d\geq 2$, there are several  big Degs.
An explicit $\Deg$, for all dimensions,  was introduced in
\cite{hdeg}. It has a recursive aspect.

\begin{Definition}\label{hdegdef} {\rm Let  $M$ be a finitely generated graded module
over the graded
 algebra $A$ and  $S$  a 
 Gorenstein graded algebra mapping onto $A$, with maximal graded ideal
${\mathfrak m}$.  Set $\dim S=r$,
 $\dim M= d$. 
The {\em homological degree}\index{homological degree}\index{hdeg,
the homological degree} of $M$\label{hdegnu}
is the integer
\begin{eqnarray} 
 \hdeg(M) &=& \deg(M) +
  \label{homologicaldegree} \\
&& \sum\limits_{i=r-d+1}^{r} {{d-1}\choose{i-r+d-1}}\cdot
 \hdeg(\mbox{\rm Ext}^i_S(M,S)).\nonumber 
\end{eqnarray}
This expression becomes more compact when $\dim M=\dim S=d>0$:
\begin{eqnarray}
 \hdeg(M) & = &\deg (M) + \label{homologicaldegree2}\\
&&
 \sum\limits_{i=1}^{d} {{d-1}\choose{i-1}}\cdot
 \hdeg(\mbox{\rm Ext}^i_S(M,S)). \nonumber
\end{eqnarray}
}\end{Definition}

We are going to recall some of the bounds afforded by a $\Deg $
function.

\begin{Theorem}[{\rm {\cite[Theorem 2.94]{icbook}}}]
\label{Degandbetti}
 For any $\Deg$ function and any finitely generated
$R$-module $M$,
\[ \beta_i(M) \leq \
\Deg(M)\cdot  \beta_i(k), \] 
where $\beta_i(\cdot)$ is the $i$th Betti number function.
\end{Theorem}

\begin{Theorem}[{\rm \cite{Nagel}}] \label{Nagel}
Let $A$ be an standard graded algebra over an infinite field and let
$M$ be a nonzero  finitely generated graded $A$-module. Then for any
$\Deg(\cdot)$ function, we have
\[ \reg(M)< \Deg(M)+ \alpha(M),\]
where $\alpha(M)$ is the maximal degree in a minimum graded
generating set of $M$, and $\reg(M)$ is its Castelnuovo-Mumford
regularity.
\end{Theorem}

\subsubsection*{Higher cohomology modules} 

Given a Gorenstein local ring
$R$, and two finitely generated $R$-modules $A$ and $B$, we look at
the problem of bounding the torsion of the modules
$\Tor_i^R(A,B)$, for $i>0$. The approach we use is
 straightforward: Consider a free presentation,
\[ 0 \rar L \lar F \lar A \rar 0,\]
pass to $L$ the given $\hdeg$ information on $A$, and use d\'ecalage
to compare $h_0(\Tor_i^R(A,B))$ to
 $h_0(\Tor_{i-1}^R(L,B))$.
This is allowed since by the  cohomology exact sequences, we have the
short exact sequences
\[ 0 \rar \Tor_1^R(A,B)\lar  L\otimes B \lar F\otimes B \lar A\otimes B  \rar 0,
\] and
\[ \Tor_i^R(A,B) \simeq \Tor_{i-1}^R(L,B), \quad i>1. \]

We relate the degrees of $A$ to those of $L$. We shall assume that
the rank of $F$ is $\nu(A)$,
$\dim A=\dim R=d$, so that if $L\neq 0$, $\dim L=d$. This gives
\[ \deg(L) = \deg(F)-\deg(A).\]

We have the two expressions for $\hdeg(A)$ and $\hdeg(L)$
(\cite[Definition 2.77]{icbook}):

\begin{eqnarray*}
 \hdeg(A) & = &\deg (A) +
 \sum_{i=1}^{d} {{d-1}\choose{i-1}}\cdot
 \hdeg(\mbox{\rm Ext}^i_R(A,R)),\\
 \hdeg(L) & = &\deg (L) +
 \sum_{i=1}^{d-1} {{d-1}\choose{i-1}}\cdot
 \hdeg(\mbox{\rm Ext}^i_R(L,R)),
\end{eqnarray*} since $\depth L> 0$. If we set $a_i=\hdeg
(\Ext_R^i(A,R))$, these formulas can be rewritten as
\begin{eqnarray*}
 \hdeg(A) & = &\deg (A) + a_1 + 
 \sum_{i=2}^{d} {{d-1}\choose{i-1}}\cdot
 a_i,\\
 \hdeg(L) & = &\deg (L) +
 \sum_{i=1}^{d-1} {{d-1}\choose{i-1}}\cdot
 a_{i+1} = \deg(L) + \sum_{i=1}^{d-1} \frac{i}{d-i}{{d-1}\choose{i}}\cdot
 a_{i+1},
\end{eqnarray*}
where we have made use of the isomorphism $\Ext_R^i(A,R)\simeq
\Ext_R^{i-1}(L,R)$, for $i>1$.

This gives
\[
 \frac{1}{d-1} \sum_{i=2}^{d}
{{d-1}\choose{i-1}}\cdot a_i\leq \hdeg(L)-\deg(L)
\leq  (d-1)\cdot \sum_{i=2}^{d}
{{d-1}\choose{i-1}}\cdot a_i.
\]

We now collect these estimations:

\begin{proposition} \label{hdegsyz} Let $R$, $A$ and $L$ be as above.
Then \[
\deg(L) \leq  \nu(A)\cdot \deg(R)-\deg(A),\]
and if $c=\hdeg(A)-\deg(A)-\hdeg(\Ext^1_R(A,R))$
\[ \frac{1}{d-1} c
\leq \hdeg(L)-\deg(L) \leq
{(d-1)} c.
\] 
\end{proposition}

\begin{Theorem} \label{h0highcoh}
Let $R$ be a Gorenstein local ring. If there is a polynomial
$\ff(x,y)$ such that for any two finitely generated $R$-modules
$A,B$, in a certain class of modules,  $h_0(A\otimes_RB)\leq \ff(\hdeg(A), \hdeg(B))$, there are also
polynomials $\ff_i(x,y)$, $i\geq 1$, of the same degree, such that 
\[ h_0(\Tor_i^R(A,B)) \leq \ff_i(\hdeg(A), \hdeg(B)).\]

\end{Theorem}

\demo Consider  a minimal free presentation of $A$,
\[ 0 \rar L \lar F \lar A \rar 0.\]
Tensoring with $B,$ we have the exact sequence
\[ 0 \rar \Tor_1^R(A,B) \rar L\otimes B \lar F\otimes B \lar A\otimes B
\rar 0,\]
hence $h_0(\Tor_1^R(A,B)) \leq h_0(L\otimes B)$. Now we make use of
Proposition~\ref{hdegsyz} to bound $h_0(L\otimes B)$ using the data
on $A$. 

The bounds for the higher Tor comes from the d'\'ecalage. \QED

\section{Dimension $1$} Suppose $R$ is a local domain of
dimension $1$. We start our discussion with the case of two ideals,
$I,J\subset R$. Consider the commutative diagram   
\[
\diagram 
0 \rto  & L  \rto\dto & F \rto \dto & IJ
\rto\dto  &0 \\
0 \rto          &  T \rto                 & I\otimes_RJ \rto  & IJ \rto  &0
\enddiagram
\]
where $F$ is a free presentation of $I\otimes_RJ$, and therefore it
has rank $\nu(I)\cdot \nu(J)$. $L$ is a torsion free  module with
$\rank(L)=\rank(F)-1$,  and therefore it can be generated by $\deg(L) =
(\rank(F)-1)\deg(R)$ elements. 

\begin{proposition} If $R$ is a   local domain of dimension one,
essentially of finite type over a field   and $I,J$ are
$R$-ideals, then
\begin{eqnarray*}
h_0(I\otimes_RJ)  &\leq & (\nu(I)\cdot \nu(J)-1)\deg (R)
\cdot \lambda(R/(I,J,K))\\
\nu(H_{\mathfrak{m}}^0(I\otimes_R J)) &\leq & (\nu(I)\cdot \nu(J)-1)\deg
(R),
\end{eqnarray*}
where $K$ is the Jacobian ideal of $R$.
\end{proposition}

We now quote in full two results \cite{Wang} that we require for the
proof.

\begin{Theorem}[{\cite[Theorem 5.3]{Wang}}] \label{Wang1} Let $R$ be
a Cohen-Macaulay local ring of dimension $d$, essentially of finite
type over a field, and let $J$ be its Jacobian ideal. Then $J\cdot
\Ext_R^{d+1}(M, \cdot) =0$ for any finitely generated $R$-module $M$,
or equivalently, $J\cdot \Ext_R^{1}(M, \cdot) =0$ for any finitely
generated maximal Cohen-Macaulay $R$-module $M$.
\end{Theorem}

\begin{proposition}[{\cite[Proposition 1.5]{Wang}}] \label{Wang2}
Let $R$ be a commutative ring, $M$ and $R$-module, and $x\in R$. If
$x\cdot \Ext_R^1(M, \cdot)=0$ then
$x\cdot \Tor_1^R(M, \cdot)=0$.
\end{proposition}

\demo The second formula arises because $L$ maps onto $T$, the torsion submodule of
$I\otimes_RJ$, which is also annihilated by $I, J$
and $K$. The last assertion
is a consequence of the fact that $T= \Tor_1^{R}(I,R/J)$, and by
Theorem~\ref{Wang1} and Proposition~\ref{Wang2}, $K$ will annihilate
it. \QED

\bigskip

A small enhancement occurs since one can replace $I$ and $J$ by
isomorphic ideals.
In other words, the last factor, $\lambda(R/(I,J,K))$, can be replaced
by  $\lambda(R/(\tau(I),\tau(J),K))$, where $\tau(I) $ and $ \tau(J)$
are their {\em trace} ideals ($\tau(I) = \mbox{\rm image
$I\otimes_R\Hom_R(I,R)\rar R$}$).

\medskip

The version for modules is similar:

\begin{proposition} If $R$ is a   local domain of dimension one,
essentially of finite type over a field   and let $A,B$ be finitely
generated  torsion free  
$R$-modules. Then
\begin{eqnarray*}
h_0(A\otimes_RB)  &\leq & (\nu(A)\cdot \nu(B)-\rank(A)\cdot \rank(B))\deg (R)
\cdot \lambda(R/K),\\
\nu(H_{\mathfrak{m}}^0(A\otimes_R B)) &\leq &
(\nu(A)\cdot \nu(B)-\rank(A)\cdot \rank(B))\deg (R),
\end{eqnarray*}
where $K$ is the Jacobian ideal of $R$.
\end{proposition}

To extend this estimation of $h_0(A\otimes_RB)$ to finitely generated torsion free
$R$-modules that takes into account annihilators
we must equip the modules--as is the case of ideals--with a
privileged embedding into free modules.

\begin{proposition} Let $R$ be a Noetherian integral domain of dimension
$1$, with finite integral closure. Let $A$ be a torsion free
$R$-module of rank $r$ with the embedding $A\rar F=R^r$. Let $I$ be
the ideal $\mbox{\rm image }(\wedge^rA \rar \wedge F = R)$. Then $I$
annihilates $F/A$.
\end{proposition}

\section{Vector bundles} Let $(R, \mathfrak{m})$ be a regular local
ring of dimension $d\geq 2$. For a vector bundle $A$ (that is, a
finitely generated $R$-module that is free on the punctured
spectrum), we establish estimates of the form
\begin{eqnarray} \label{vbtorsion}
h_0(A \otimes_R B) &\leq & c(R) \cdot \hdeg(A) \cdot  \hdeg(B), 
\end{eqnarray}
where $c(R)$ is a constant depending on $R$.

\medskip

From the general observations above, we may assume that $\depth A$
and $\depth B$ are positive.

\medskip

We make some reductions beginning with the following. Since $A$ is
torsion free, consider the natural exact sequence
\[ 0 \rar A \lar A^{**} \lar C \rar 0.
\] Note that
$C$ is a module of finite support, and that $A^{**}$ is a vector
bundle. Furthermore, since $A^{**}$ has depth at least $2$, a direct
calculation yields
\[ \hdeg(A) = \hdeg(A^{**}) + \hdeg(C).\]

Tensoring the exact sequence by $B$, gives the exact complex
\[ \Tor_1^R(C,B) \lar A\otimes_RB \lar 
A^{**} \otimes_RB,\]
from which we obtain
\[ h_0(A\otimes_RB) \leq  h_0(A^{**}\otimes_RB) +
\lambda(\Tor_1^R(C,B)).
\]
As $C$ is a module of length $\hdeg(C)$,
\[ \lambda(\Tor_1^R(C,B))\leq \beta_1(k)\cdot \hdeg(C)\cdot \nu(B).
\] We recall that $\nu(B)\leq \hdeg(B)$. 

\medskip

\bigskip

Let $(R, \mathfrak{m})$ be a Gorenstein local ring of dimension $d>0$
and let $A$ be a finitely generated $R$-module that is free on the
punctured spectrum, and has finite projective dimension.
 We seek to estimate $h_0(A\otimes B)$ for various $R$-modules $B$.

\begin{Theorem}\label{vbbetahi} $B$ be a module of projective dimension $< d$ and let
$A$ be a module that is
free on the punctured spectrum.
Then
\[ h_0(A\otimes B) \leq \sum_{i=0}^{d-1} \beta_i(B)\cdot h_i(A).\]
\end{Theorem}

\demo Let $B$ be an $R$-module with $\depth B>0$ 
with the minimal free resolution 
\[ 0\rar F_{d-1} \lar F_{d-2} \lar \cdots \lar F_1 \lar F_0 \lar B
\rar 0. \]
Tensoring by $A$ gives a complex
\[ 0\rar F_{d-1}\otimes A \lar F_{d-2}\otimes A  \lar \cdots \lar
F_1\otimes  A\lar F_0\otimes A \lar B\otimes A \rar 0, \]
whose homology $H_i=\Tor_i(B,A)$, $d>i>0$, has finite support. Denoting by $B_i$
and $Z_i$ its modules of boundaries and cycles, 
we have several exact sequences that start at $Z_{d-1}=0$
\[ 0 \rar B_0 \lar F_0\otimes A \lar B\otimes A \rar 0,\]
\[ 0 \rar B_i \lar Z_i \lar H_i \rar 0,\]
\[ 0 \rar Z_{i} \lar F_{i}\otimes A \lar B_{i-1} \rar 0.\]
Taking local cohomology, we obtain the following acyclic complexes of
modules of finite length
\[ 
 H^0_{m}(F_0\otimes A) \lar H^0_{m}(B\otimes A) \lar
H^1_{m}(B_0), \]

\[  H^i_{m}(F_{i}\otimes A) \lar H^i_{m}(B_{i-1})
\lar H^{i+1}_{m}(Z_{i}), \quad d-1>i\geq 1. \]

Now we collect the inequalities of length, starting with
\[ h_0(B\otimes A)\leq \beta_0(B)\cdot h_0(A) + h_1(B_0)
\]
and
\[ 
 h_i(B_{i-1})\leq \beta_i(B)\cdot h_i(A) + h_{i+1}(Z_{i})
\leq \beta_i(B)\cdot h_i(A) + h_{i+1}(B_{i}),
\]
where we replace the rank of the modules $F_i$ by $\beta_i(B)$. \QED


\begin{Theorem} \label{h0vb} Let $R$ be a regular local ring of
dimension $d$. If $A$ is a finitely generated module free on the punctured
spectrum, then for any finitely $R$-module $B$
\begin{eqnarray*}
h_0(A\otimes B) \leq d\cdot \hdeg(A)\cdot \hdeg(B).
\end{eqnarray*}
\end{Theorem}

\demo
Let us rewrite the inequality in Theorem~\ref{vbbetahi}
 in case $R$ is a regular local ring.
Since by local duality $h_i(A)= \lambda(\Ext_R^{d-i}(A,R))$ and
$\beta_i(B)\leq \beta_i(k)\cdot \hdeg(B)$, we have
\begin{eqnarray*} h_0(A\otimes B) &\leq  & \hdeg(B )\cdot \sum_{i=0}^{d-1}{d\choose
i}\lambda(\Ext_R^{d-i}(A,R))\\
&\leq & d\cdot \hdeg(B)\sum_{i=1}^{d} 
 {{d-1}\choose{i-1}}
\lambda(\Ext_R^{d-i}(A,R))\\
& \leq & d\cdot \hdeg(A)\cdot \hdeg(B).
\end{eqnarray*}

Finally, to deal with a general module $B$,  it suffices 
to add to $\hdeg(B)$ the correction $h_0(B)$ as given in
Proposition~\ref{AmodA0}. \QED

\begin{Example} \label{bv} {\rm Let $(R, \mathfrak{m})$ be a Gorenstein local
ring of dimension $d\geq 1$ and let $A$ be a module with a
presentation
\[ 0 \rar R^n \stackrel{\varphi}{\lar} R^{n+d-1} \lar A \rar 0,\]
where the ideal $I_n(\varphi)$ is $\mathfrak{m}$-primary. $A$ is a
vector bundle of projective dimension $1$.  
According to a well-known length formula (\cite{BV86}), $\lambda(\Ext_R^1(A,R))= \lambda
(R/I_n(\varphi))$, from which it follows that 
\begin{eqnarray*}
\hdeg(A) &=& (d-1)\deg(A) + \lambda(R/I_n(\varphi)), \quad \mbox{\rm
and therefore}\\
h_0(A\otimes A) &\leq & d((d-1)\deg(A) + \lambda(R/I_n(\varphi)))^2.
\end{eqnarray*}
}\end{Example}


\section{Dimension $2$}
Let $(R, \mathfrak{m})$ be a regular local ring of dimension $2$ (or
a polynomial ring $k[x,y]$ over the field $k$). For two $R$-modules
$A$ and $B$ we are going to study
$ h_0(A\otimes_RB) = \lambda(H^0_{\mathfrak{m}}(A\otimes_RB))$
through a series of reductions on $A$, $B$ and $R$. Now we examine
how the presence of torsion affects the analysis.

Eventually the
problem will settle on the consideration of a special class of
one-dimensional rings.

\medskip

We are already familiar with the stripping away from $A$ and $B$ of their
submodules of finite support, so we may assume that these modules
have depth $\geq 1$. Let $A$ be a module of dimension $2$,
and denote by $A_0$ the torsion submodule of
$A$. Consider the natural exact sequence
\[0 \rar A_0 \lar A \lar A'\rar 0,\]
with $A'$ torsion free. If $A_0\neq 0$, it is a Cohen-Macaulay module
of dimension $1$. We have the exact sequence
\[ 0 \rar \Ext_R^1(A',R) \lar \Ext_R^1(A,R) \lar \Ext_R^1(A_0,R)
\rar 0,
\] that yields
 \begin{eqnarray*}
\deg A &=& \deg A' \\
\hdeg(\Ext_R^1(A,R)) &=& 
\hdeg(\Ext_R^1(A',R)) +\deg A_0,
\end{eqnarray*} 
in particular
\begin{eqnarray*}
\hdeg(A) &=& \deg A_0 + \hdeg(A').
\end{eqnarray*}

As a consequence, $\hdeg(A')$ and $\hdeg(A_0)$ are bounded  in
terms of $\hdeg(A)$. Now we tensor the sequence by $B$ to get the
complex
\[ \Tor_1^R(A',B) \lar A_0\otimes_R B \lar 
 A\otimes_R B \lar  A\otimes_R B \rar 0. \]
If we denote by $L$ the image of $A_0\otimes_R B$ in $A\otimes_R B$,
since $\Tor_1^R(A',B)$ is a module of finite length, we
have
\begin{eqnarray*} h_0(A\otimes_RB) &\leq & h_0(A'\otimes_RB) + h_0(L),\\
h_0(L) & \leq & h_0(A_0\otimes_RB). 
\end{eqnarray*}

If we apply a similar reduction to $B$ and combine, we get
\begin{eqnarray*} h_0(A\otimes_RB) &\leq h_0(A'\otimes_RB') + 
h_0(A_0\otimes_RB') + 
h_0(A'\otimes_RB_0)
+  h_0(A_0\otimes_RB_0). 
\end{eqnarray*}

A term like $h_0(A'\otimes B_0)$ is easy to estimate since $A'$ is a
vector bundle of low dimension. For convenience, let 
\[ 0 \rar F_1 \lar F_0 \lar A' \rar 0\] be a minimal resolution of
$A'$. Tensoring with $B_0$, we get the exact sequence
\[ 0 \rar \Tor_1^R(A',B_0) \lar F_1\otimes_RB_0
\lar F_0\otimes_RB_0
\lar A'\otimes_RB_0 \rar 0.
\] Since 
$\Tor_1^R(A',B_0)$
has finite support and $F_1\otimes_RB_0$ has positive depth,
$\Tor_1^R(A',B_0)=0$. To compute $h_0(A'\otimes_RB_0)$, consider 
a minimal resolution of $B_0$,
\[ 0 \rar G \lar G \lar B_0 \rar 0,\]
and
 the exact sequence
\[ 0 \rar A'\otimes_R G \lar A' \otimes_R G \lar A'\otimes_RB_0 \rar 0.
\] From the cohomology exact sequence, we have the surjection
\[ \Ext_R^1(A'\otimes G,R) \lar \Ext_R^2(A'\otimes_R B_0),\]
and therefore, since $\Ext_R^1(A,R)$ is a module of finite support, 
\[ h_0(A'\otimes_RB_0) \leq \nu(G)\cdot \lambda(\Ext_R^1(A',R)).\]
This shows that
\[h_0(A'\otimes_RB_0) \leq \nu(G)\cdot (\hdeg(A')-\deg A') <
\hdeg(A')\cdot \hdeg(B_0)
.\]

The reductions thus far lead us to assume that $A$ and $B$ are
$R$-modules of positive depth and dimension $1$. Let 
\[ 0 \rar F \stackrel{\varphi}{\lar} F \lar A \rar 0\]
be a minimal free resolution of $A$. By a standard calculation, 
\[ \deg A = \deg(R/\det(\varphi)).\]
Since $\det (\varphi)$ annihilates $A$, we could view $A\otimes_RB $
as a module of over $R/(\det(\varphi \circ \psi))$ where $\psi$ is
the corresponding matrix in the presentation of $B$.

\medskip

To avoid dealing with two matrices, replacing $A$ by $A\oplus B$, we
may consider $h_0(A\otimes_RA)$, but still denote by $\varphi$ the
presentation matrix (instead of $\varphi\oplus \psi$), and set
$S=R/(\det(\varphi))$; note that $\deg S = \deg A$.

\begin{Example}{\rm 
We consider a cautionary family of examples to show that other
numerical readings must be incorporated into the estimates for
$h_0(A\otimes_RA)$.

Let $A$ be a module generated by two elements, with a free resolution
\[ 0 \rar F \stackrel{\varphi}{\lar} F \lar A \rar 0.\]
Suppose $k$ is a field of characteristic $\neq 2$. To calculate
$h_0(A\otimes_RA)$, we make use of the decomposition 
\[ A \otimes_RA = S_2(A) \oplus \wedge^2A.\] 
Given a matrix representation,
\[ \varphi = \left[ \begin{array}{ll}
a_{11} & a_{12} \\
a_{21} & a_{22} \\
\end{array} \right],
\]
one has 
\[\wedge^2 A \simeq R/I_1(\varphi)=R/(a_{11}, a_{12}, a_{21},
a_{22}).\]

The symmetric square of $A$, $S_2(A)$, has a resolution
\[ 0 \rar R \stackrel{\phi}{\lar} F\otimes_RF \stackrel{\psi}{\lar}
S_2(F),\]
where 
\begin{eqnarray*}
\psi(u\otimes v) &=& u \cdot  \varphi(v) + v\cdot  \varphi(u) \\
\phi (u\wedge v) & = & \varphi'(u)\otimes v-\varphi'(v)\otimes u \\
\end{eqnarray*}
where $\varphi'$ is the matrix obtained from $\varphi$ by dividing
out its entries by their gcd $a$, $\varphi=a\cdot \varphi'$.

\medskip

A straightforward calculation will give
\[ \Ext_R^2(S_2(A), R) = R/I_1(\varphi').\]
This shows that
\begin{eqnarray*}
 h_0(A\otimes_RA) & =& h_0(R/I_1(\varphi))+h_0(R/I_1(\varphi'))\\
&=& 
h_0(aR/aI_1(\varphi'))+h_0(R/I_1(\varphi'))= 
2\cdot \lambda(R/I_1(\varphi')).
\end{eqnarray*}

Thus  the matrix 
\[ \varphi = \left[ \begin{array}{ll}
x & y^n \\
0 & x \\
\end{array} \right],
\]
will define a module $A$, with $\deg(A)=2$, but
$h_0(A\otimes_RA)=2n$. This means that we must take into account the
degrees of the
entries of $\varphi$ itself.

}\end{Example}

\begin{Example}\label{degabc}{\rm Let $R$ be a Cohen-Macaulay local
ring and let $\{ x_1, \ldots, x_n\}$, $n\geq 2$,  be a set of elements such that
any pair forms a regular sequence. Set
\begin{eqnarray*}
 \xx &= & x_1\cdots x_n,\\
z_i & = & \xx/x_i, \quad i=1\ldots n.
\end{eqnarray*}
We claim that 
\begin{eqnarray} \label{degxx}
 \deg(R/(z_1, \ldots, z_n)) \leq \frac{1}{2} ((\deg(R/(\xx)))^2-n).
\end{eqnarray}
We argue by induction on $n$, the formula being clear for $n=2$.

Consider the exact sequence
\[ 0 \rar (x_1, z_2, \ldots, z_n )/(z_1, z_2, \ldots, z_n) \lar
R/(z_1, z_1) \lar R/(x_1,z_1) \rar 0.\]
Since 
\[ (x_1,z_1)/(z_1, \ldots, z_n)\simeq R/(z_2', \ldots, z_n'),
\] where $z_i'$, $i\geq 2$, denotes the products from elements in the set $\{x_2,
\ldots, x_n\}$ using the formation rule of the $z_i$.

Adding the multiplicities of the modules of the same dimension, we
have
\[ \deg(R/(z_1, \ldots, z_n)) = \deg(R/(z_1,\ldots, z_n)) + \deg(R/(z_2',
\ldots, z_n')). \]
As
 \[\deg(R/(x_1, z_1))= \deg(R/(x_1))\cdot \deg(R/(z_1))= 
\deg(R/(x_1))\cdot \sum_{j\geq 2}\deg(R/(x_j)) ,\] and by
induction
\[ \deg(R/(z_2', \ldots, z_n'))= \sum_{2\leq i< j\leq
n}\deg(R/(x_i))\cdot \deg(R/(x_j)),\]
we have
\[ \deg(R/(z_1, \ldots, z_n) = \sum_{1\leq i< j \leq n}
\deg(R/(x_i))\cdot \deg(R/(x_j)).
\]

The rest of the calculation is clear.
There are similar formulas in case every subset of $k$ elements of
$\{x_1, \ldots, x_n\}$ forms a regular sequence.

}\end{Example}


Now we return to the modules with a presentation
\[ 0 \rar F \stackrel{\varphi}{\lar} F \lar A \rar 0,\]
and write $\det (\varphi)=\xx=x_1\cdots x_n$.
 Setting $z_i= \xx/x_i$, consider the exact
sequence
\[ 0 \rar R/(\xx) \lar R/(z_1) \oplus \cdots \oplus R/(z_n) \lar C
\rar 0,\]
induced by the mapping $1 \mapsto (z_1, \ldots, z_n)$. $C$ is a
module of finite length, and making use of duality and the inequality
(\ref{degxx}), 
\[ \lambda(C) \leq \frac{1}{2} ((\deg(R/(\xx))^2-n).\]

Tensoring this sequence by $A$, gives
\[ \Tor_1^R(A,C) \lar A \lar A_1 \oplus \cdots \oplus A_n \lar
A\otimes_RC \rar 0,\]
and since $\depth A>0$, we have the exact sequence
\[ 0 \rar A\lar A_1 \oplus \cdots \oplus A_n \lar A\otimes_RC  \rar 0,\]
where $A_i = A/x_iA$ and $\lambda(A\otimes C)\leq \nu(A) \lambda(C)$.
These relations give that 
\begin{eqnarray*}
\deg A &=& \sum_{i=1}^n \deg A_i \\
\lambda(A\otimes_R C) &\geq & \sum_{i=1}^n h_0( A_i). 
\end{eqnarray*}
These inequalities show that we are still tracking the $\hdeg(A_i)$ in
terms of $\deg A$.

Tensoring the last  exact sequence by $A$, we obtain the exact
complex
\[ \Tor_1^R(A, A\otimes_RC) \lar A\otimes_RA \lar A_1\otimes_R A_1
\oplus \cdots \oplus A_n\otimes_RA_n, \]
from which we have
\begin{eqnarray*} h_0(A\otimes_RA) &\leq & 
\sum_{i=1}^n h_0(A_i\otimes_RA_i) + \lambda
(\Tor_1^R(A, A\otimes_RC)) \\
&\leq & \sum_{i=1}^n h_0(A_i\otimes_RA_i) + \beta_1(A)\cdot \nu(A)
\cdot \lambda(C). \\
\end{eqnarray*}

Let us sum up these reductions as follows:

\begin{proposition} Let $R$ be a two-dimensional regular local ring
and let $A$ be a Cohen-Macaulay $R$-module of dimension one. Then
\[ h_0(A\otimes_RA) \leq 3 \cdot  \hdeg(A)^4,\]
provided 
\[ h_0(A\otimes_RA) \leq 2 \cdot  \hdeg(A)^4\]
 whenever $\ann A$ is a primary ideal.
\end{proposition}

\demo Note that $\beta_1(A) \leq \beta_1(k) \cdot \hdeg(A)$,  $\nu(A)\leq
\hdeg(A)$, and $\lambda(C)< \frac{1}{2} \hdeg(A)^2$.
\QED

\medskip

\section{Dimension $3$} The technique of 
Theorem~\ref{h0vb} can be used to deal with torsionfree modules of
dimension three.  

\begin{Theorem} \label{h0dim3}
Let $R$ be a regular local ring of dimension $3$,
and let   $A$ and $B$ be torsionfree $R$--modules.
Then
\begin{eqnarray}
h_0(A\otimes B)& < & 4\cdot \hdeg(A)\cdot \hdeg(B).
\end{eqnarray}
\end{Theorem}

\demo
 Consider the
natural exact sequence
\[ 0 \rar A \lar A^{**} \lar C \rar 0.\]
A straightforward calculation will show that
\begin{eqnarray} \label{h0eq0}
\hdeg(A) &=& \hdeg(A^{**}) + \hdeg(C).
\end{eqnarray}
Note that $A^{**}$ is a vector bundle of projective dimension at most
$1$ by the Auslander-Buchsbaum equality (\cite[Theorem 1.3.3]{BH}), and $C$ is a
module of dimension at most $1$. Tensoring by
the torsionfree $R$-module $B$, we have
the exact sequence
\[ \Tor_1^R(A^{**}, B) \lar \Tor_1(C,B) \lar A\otimes B \lar
A^{**}\otimes B \lar C\otimes B \rar 0,
\] where $\Tor_1^R(A^{**}, B)=0$, since $\mbox{\rm proj dim
}A^{**}\leq 1$ and $B$ is torsionfree.

\medskip

From the exact sequence, we have
\begin{eqnarray} \label{h0eq1}
h_0(A\otimes B) &\leq & h_0(A^{**}\otimes B) + h_0(\Tor_1^R(C,B)).
\end{eqnarray}
Because $A^{**}$ is a vector bundle, by Theorem~\ref{h0vb},
\begin{eqnarray} \label{h0eq2}
h_0(A^{**}\otimes B) &\leq & 3\cdot
\hdeg(A^{**})\cdot \hdeg(B).
\end{eqnarray}

For the module $\Tor_1^R(C,B)$, from a minimal free presentation of
$B$,
\[ 0 \rar L \lar F \lar B\rar 0,\]
we have an embedding $\Tor_1^R(C,B) \rar C\otimes L$, and therefore
\begin{eqnarray} \label{h0eq3}
h_0(\Tor_1^R(C,B)) &\leq & h_0(C\otimes L) \leq 3 \cdot \hdeg(L)\cdot
\hdeg(C),
\end{eqnarray}
because $L$ is a vector bundle. In turn
\begin{eqnarray*}
\deg(L) &=& \beta_1(B)-\deg(B)\leq
\beta_1(R/\mathfrak{m})\cdot \hdeg(B)-\deg(B) \\
&=& 3\cdot \hdeg(B)-\deg(B)
\end{eqnarray*}
by Theorem~\ref{Degandbetti},
and since $\Ext_R^1(L,R)= \Ext_R^2(B,R)$, 
 \begin{eqnarray*}
\hdeg(L) = \deg(L)+ \hdeg(\Ext_R^1(L,R)) & < & 4\cdot \hdeg(B).
\end{eqnarray*}

Finally we collect (\ref{h0eq2}) and (\ref{h0eq3}) into
(\ref{h0eq1}), along with (\ref{h0eq0}),
\begin{eqnarray*}
h_0(A\otimes B) &< & 3\cdot \hdeg(A^{**})\cdot \hdeg(B) + 4\cdot
\hdeg(B)\cdot \hdeg(C) \\
 &< & 4\cdot \hdeg(A) \cdot \hdeg(B),
\end{eqnarray*}
as asserted.
\QED

\section{Graded modules} 
We give a rough (high degree) estimate for the case of
graded modules over $R=k[x,y]$. We may assume that $A$ is not a
cyclic module. Furthermore, we shall assume that $A$ is
equi-generated.

\medskip

We briefly describe the behavior of $\reg(\cdot)$ with regard to some
exact sequences.

\begin{proposition}\label{addiofreg} Let $R$ be a standard graded algebra, and let
\[ 0 \rar A \lar B \lar C \rar 0\] be an exact sequence of finitely
generated graded $R$-modules, then
\[ \reg(B)\leq \reg(A)+ \reg(C),
\] 
Similarly, 
\[ \reg(A)\leq \reg(B)+ \reg(C), \quad \mbox{\rm and  }\quad
\reg(C)\leq \reg(A)+ \reg(B).\]
\end{proposition}

Let us apply it to the graded $k[x,y]$-module $A$ of depth $> 0$. In
the exact sequence
\[ 0 \rar A_0 \lar A \lar A'\rar 0 \]
we already remarked that $\hdeg(A)=\hdeg(A') + \deg A_0$. It is also
the case that if $A$ is generated by elements of degree $\leq
\alpha(A)$, then by the proposition above and Theorem~\ref{Nagel},
$\reg(A_0)< \hdeg(A) + \alpha(A)$. Actually, since $A_0$ is
Cohen-Macaulay,  a direct calculation will show that $\reg(A_0)\leq
\reg(A)$.

\medskip

We may assume that $A$ is a one-dimensional graded $R$-module with a
minimal resolution
\[ 0 \rar F \stackrel{\varphi}{\lar} F \lar A \rar 0.\]
A presentation of $A\otimes_RA$ is given by
\[ F\otimes_RF \oplus F\otimes_RF \stackrel{\psi}{\lar} F\otimes_RF,  \]
where 
$\psi = \varphi\otimes I - I\otimes \varphi$. The kernel of $\psi$ contains the
image of 
\[ \phi: F\otimes_RF \lar F\otimes_RF \oplus F\otimes_RF, \quad \phi=
I \otimes \varphi \oplus \varphi \otimes I.\]

Since $R=k[x,y]$, $L=\ker(\psi)$ is a free $R$-module of rank $r^2$,
$r=\rank(F)$. Because $\phi$ is injective, its image $L_0$ is a free
$R$-submodule of $\FF=F\otimes_R F \oplus F\otimes_R F$ of the same rank as
$L$, $L_0\subset L$, $\phi': F\otimes F\rar \FF$.
 It follows that the degrees of the entries of
$\phi'$ cannot be higher than those of $\phi$.

To estimate  $h_0(A\otimes_RA)$, note that
$\Ext_R^2(A\otimes_RA,R)$ is the cokernel of map  $\phi$. This is a
module generated by $r^2$ elements, annihilated by the maximal minors
of $\phi$.
We already have that $\ff=\det(\varphi)$ annihilates $A$. Now we look for
an element  $\hh$ in the ideal of maximal minors of $\phi'$
 so that $(\ff, \hh)$ has finite colength, and as a
consequence we would have
\[ h_0(A\otimes_RA) \leq r^2\cdot \lambda(R/(\ff,\hh))=
r^2\cdot \deg(R/(\ff))\deg(R/(\hh)).\]

The entries of $\varphi$ have degree $\leq \deg(A)-1$, so the minors
$\hh$ of $\phi$ have degree 
\[\deg \hh \leq r^2\cdot (\deg(A)-1). 
\]

\begin{proposition} \label{grdim2} If $A$ is a graded $k[x,y]$-module of dimension $1$,
equigenerated in degree $0$,
then
\[ h_0(A\otimes_RA) \leq r^4\cdot \deg(A)(\deg(A)-1)< \deg(A)^6.\]
\end{proposition}

\section{Some open questions}

There are numerous open issues regarding the torsion in tensor
products that are not discussed here.  We raise a few related to the
discussion above.

\begin{enumerate}

\item How good are some of the estimates for $h_0(A\otimes A)$
compared to actual values, for instance of Example~\ref{bv}?

\item In \cite{DV2}, Theorem~\ref{Wang1} is used to extend the HomAB
version of Theorem~\ref{h0vb} from the regular to the isolated
singularity case. Is there a similar extension to $h_0(A\otimes B)$? 

\item How to derive more general estimates in dimension $3$,
particularly of graded modules with torsion? 

\item Let $R$ be a Noetherian local domain and $A$ a finitely
generated torsionfree $R$-module. Is there an integer $e=e(R)$ 
guaranteeing that if $M$ is not $R$-free, then the tensor power
$M^{\otimes {e}}$ has nontrivial torsion? The motivation is 
 a result of Auslander
(\cite{Aus61}, see also \cite{Licht66}) that asserts that $e=\dim R$ works
for all regular local rings.
For instance, if $R$ is a one-dimensional domain, will $e=2$ work?
A more realistic question is, if $R$ is a Cohen-Macaulay local domain
of dimension $d$ and multiplicity $\mu$, will
\[ e= d+\mu-1\]
suffice? Note that if we make no attempt to determine uniform bounds
for $e$, if $\bbq\subset R$, then for a module $M$ of rank $r$ and
minimal number of generators $n$, then the embedding
\[ 0\neq \wedge^n M \hookrightarrow M^{\otimes n}\]
shows the existence of a test power {\em for} $M$.

\end{enumerate}

\end{document}